\documentclass[12pt,a4paper]{article}
\usepackage{amsfonts,amsthm,latexsym,amsmath,amssymb,amscd,epsfig,psfrag,enumerate}
\usepackage{graphics,graphicx, bezier, float, color}
\usepackage{graphicx}
\usepackage{caption}

\newtheorem{theo+}           {Theorem}
\newtheorem{prop+}           {Proposition}
\newtheorem{coro+}           {Corollary}
\newtheorem{lemm+}           {Lemma}

\newtheorem{ques+}{Question}
\theoremstyle{definition}

\newtheorem{not+}            {Notation}

\newtheorem{Ex}              {Example}[section]

\newtheorem*{ack}            {Acknowledgements}

\newtheorem{rema+}           {Remark}

\usepackage[left=16mm, right=16mm, top=20mm, bottom=20mm] {geometry}
\newenvironment{theorem}{\begin{theo+}}{\end{theo+}}
\newenvironment{proposition}{\begin{prop+}}{\end{prop+}}
\newenvironment{corollary}{\begin{coro+}}{\end{coro+}}
\newenvironment{lemma}{\begin{lemm+}}{\end{lemm+}}
\newenvironment{remark}{\begin{rema+}}{\end{rema+}}

\newcommand{\bq}{\begin{eqnarray}}
\newcommand{\eq}{\end{eqnarray}}
\newcommand{\beq}{\begin{eqnarray*}}
\newcommand{\eeq}{\end{eqnarray*}}

\newtheorem{definition}{\bf Definition}[section]

\date{}

\begin{document}

\begin{center}

{\bf \Large Characterizations of regular modules}

 \vspace*{0.5cm}

 Philly Ivan Kimuli\footnote{Department of Mathematics, Muni University, P.O. Box 725, Arua, Uganda}  and David Ssevviiri\footnote{Corresponding author}
\footnote{Department of Mathematics, Makerere University, P.O. Box 7062, Kampala, Uganda}\\
E-mail addresses:  pi.kimuli@muni.ac.ug;  david.ssevviiri@mak.ac.ug
\end{center}

\begin{abstract}
\noindent Different and distinct notions of  regularity for modules exist in the literature. When these notions are restricted to commutative rings, they all coincide  with the well-known  von-Neumann regularity for rings.  We give new characterizations of these distinct  notions for modules  in terms of both (weakly-)morphic modules and reduced modules. Furthermore, module theoretic settings are established where these in general distinct notions turn out to be indistinguishable.
\end{abstract}



\section{Introduction}
\noindent  Let $R$  be an associative and unital ring that is not necessarily commutative and $M$ be a right $R$-module.   We call   $R$  ({\it unit-}){\it regular} if for each $a\in R$ there exists a (unit) $y\in R$ such that $a = aya$.
It is {\it strongly regular} if for each $a\in R$ there exists an element $y\in R$ such that $a=a^2y$, or equivalently if it is regular and  idempotents are central.
In the literature, there are different characterizations of a regular ring  which are  distinct for modules.
For instance,  see \cite[pg. 237]{ware1971endomorphism} and \cite[Exercises 15 (13)]{anderson1992rings},  a ring is regular $\Leftrightarrow$  every right (left) cyclic ideal is a direct summand $\Leftrightarrow$ every finitely generated right (left) ideal is a direct summand.
$R$ is strongly regular $\Leftrightarrow$ it is regular and reduced $\Leftrightarrow$ every right (left) cyclic ideal is generated by a central idempotent $\Leftrightarrow$  it is regular and  $Ra\subseteq aR$ for every $a\in R \Leftrightarrow  aR=a^2R$ for each $a\in R$.
Where $R$ is commutative, it is regular $\Leftrightarrow$ it is strongly regular.\\

\noindent Following the  (von-Neumann) regularity characterizations for rings, different authors have come up with
 different  definitions for the notion of ``regularity'' for  modules.   We outline some of them below (see also Definition~\ref{trees}, \cite[Definition 2.3]{ware1971endomorphism} and \cite{zelmanowitz1972regular}):

\begin{definition}~\label{thatis}    An  $R$-module $M$ is said to be
\begin{enumerate}[(a)]
\item {\it endoregular} \cite{lee2013modules, ware1971endomorphism}  if $\varphi(M)$ and  $ker(\varphi)$ are direct summands of $M$ for every endomorphism $\varphi$ of $M$;
\item {\it Abelian endoregular} \cite{lee2013modules}  if  $\text{End}_R(M)$ is a strongly regular ring;
\item {\it F-regular} \cite{Fieldhouse}  if for every submodule $N$ of  $M$,  the sequence $0\to N\bigotimes E\to M\bigotimes E$  is exact for each $R$-module $E$;
\item  {\it strongly F-regular} \cite{ramamurthi1973finitely}  if every finitely generated submodule of $M$ is a direct summand of $M$.\\
(In the bullets (e), (f) and (g) below, $R$ is commutative.)
\item {\it JT-regular} \cite{jayaram2018neumann} if for each  $m\in M,mR=Ma=Ma^2$ for some $a\in R$;
\item {\it weakly JT-regular} \cite{anderson2019module} if  $Ma=Ma^2$ for each $a\in R$;
\item {\it weakly-endoregular} \cite{anderson2021endoregular} if $Ma$ and  $l_M(a)$ are direct summands of $M$ for each $a\in R$.
\end{enumerate}
 \end{definition}

\noindent  An $R$-module $M$ is {\it reduced} \cite{lee2004reduced} if  whenever $a\in R$ and $m\in M$ satisfy $ma^2=0$, then $mRa =0$.           Reduced modules are a generalisation of reduced rings.  Recall that a ring  is said to be {\it reduced} \cite{marks2003taxonomy} if it has no non-zero nilpotent elements.  Thus $R$ is a   reduced ring if and only if $R$ is a reduced $R$-module.\\

\noindent   We call  $M$ a {\it morphic module} if every endomorphism $\varphi$ of $M$ has a cokernel which is isomorphic to its kernel, i.e., if for every
endomorphism $\varphi$ of $M, M/\varphi(M)\cong \ker(\varphi)$ as $R$-modules.  Note that the property  $M/\varphi(M) \cong\ker(\varphi)$  is  the dual of the  First Isomorphism Theorem for the module endomorphism $\varphi$. This notion has been widely studied, see for instance \cite{nicholson2005survey}.   Recently for commutative rings,   $M$ is called {\it weakly-morphic} in \cite{kimuli2021} if $M/Ma\cong l_M(a)$  as $R$-modules for each $a\in R$, i.e., if every endomorphism $\varphi_a$ of $M$ given by  right multiplication by $a\in R$ is morphic.  It turns out that a (commutative) ring $R$ is right (and left) morphic if and only if the $R$-module $R$ is a weakly-morphic module. Morphic rings have been studied in \cite{bamunoba2020morphic, nicholson2004rings, nicholson2004principal}.\\

\noindent The relationship between morphic, reduced and regular rings  has been extensively investigated in the literature, dating back
to when Ehrlich \cite{ehrlich1976units} proved that a ring is right morphic and regular if and only if it is unit-regular.
Since then, the study of the morphic property in rings has  flourished due to the way morphic rings connect with reduced rings to provide  conditions related to  regular rings.  Recall that a  reduced ring need not be regular or (right) morphic in general.   For example,  the ring of integers  $\Bbb Z$  is  reduced.  However, it  is neither (right) morphic nor  regular.   In general,  we have the following relations about rings:
strongly regular (i.e., regular with central idempotents) $\Leftrightarrow$ reduced and (right) morphic  $\Rightarrow$ unit-regular $\Leftrightarrow$ (right) morphic and regular $\Rightarrow$ regular.
Indeed  by  \cite[Proposition 4.13]{bamunoba2020morphic}, strongly regular (i.e.,  regular with central idempotents) rings  coincide with reduced and (right) morphic rings, and these are  unit-regular.   By \cite[Theorem 1]{ehrlich1976units},   unit-regular rings are exactly the (right) morphic and regular rings.   The remaining implications do not reverse in general.   The ring   $R:=\Bbb R\bigoplus \Bbb R$ (where $\Bbb R$ is the ring of real numbers) is unit-regular and hence morphic.   But $R$ has some nontrivial idempotents which are not central (therefore, not reduced) by  \cite[Corollary to Theorem 1]{ehrlich1976units} and \cite[Example 2.25]{lee2013modules}.  The  ring $\bigoplus_{n=1}^\infty\Bbb R$   is  regular but not right morphic (therefore, not unit-regular) and the ring $\Bbb Z/4\Bbb Z$ is morphic but not  regular.\\

\noindent This paper  gives new  characterizations of   regular modules given in Definition~\ref{thatis} in terms of (weakly-)morphic and reduced (sub-)modules.
We prove that a module is   weakly-morphic and reduced  if and only if it is weakly-endoregular (Theorem~\ref{35}); the class of Abelian endoregular modules coincides with that of morphic modules  with reduced   rings of endomorphisms (Theorem~\ref{31}); if a module $M$ is  strongly F-regular, then each of its sub-module  is invariant under every endomorphism of $M$ if and only if $M$ is a  morphic module with a reduced  ring of endomorphisms  (Theorem~\ref{3}).
A module is F-regular  if and only if each of its (cyclic) sub-modules is a weakly-morphic and reduced module (Theorem~\ref{230}).
Conditions for which one still gets coincidence of different notions of regularity  in the module theoretic setting are established.  For instance, in the subcategory of finitely generated modules, the following coincide:   weakly-morphic and reduced $\Leftrightarrow$ F-regular $\Leftrightarrow$ weakly-endoregular $\Leftrightarrow$ weakly JT-regular (Theorem~\ref{349}).\\

\noindent {\bf Notation and conventions}. Throughout this paper, all rings $R$ will be associative and unital but not necessarily commutative, $M$ is a unitary right $R$-module and $S$ denotes $\text{End}_R(M)$, the  ring of endomorphisms of $M$.   Therefore, in this  case $M$ can be viewed as a left $S$-right $R$-bimodule.      By $\mathbb{Z}, \mathbb{Q}$ and $\Bbb R$ we denote the ring of integers,  rational numbers and real numbers respectively.    For $\varphi\in S, \ker(\varphi)$ and $\text{Im}(\varphi)$ denote the kernel and image of $\varphi$ respectively.    The notation $N \subseteq M$(resp., $N\subseteq^{\bigoplus} M$) means that $N$ is a  submodule (resp., a direct summand) of $M$.
We also define $r_M(I): = \{m\in M: I(m)= 0\}, l_S(I):=\{\varphi\in S:\varphi I=0\}, r_S(I):=\{\varphi\in S:I\varphi=0\}$ for a nonempty subset $I$ of $S;~~r_R(N):=\{a\in R:Na=0\}, l_S(N):=\{\varphi\in S:\varphi(N)=0\}$  for $N \subseteq M$ and $l_M(A):=\{m\in M:mA=0\}$  for  $A\subseteq R$. Note that $r_M(\varphi):=\ker(\varphi)$ for $0\neq\varphi\in S$ and  $\text{Ann}_R(M):=r_R(M)$, the (right) annihilator of $M$.  For any $a\in R$, the principal ideal generated by $a$ is denoted by $(a)$.\\

\noindent The following definitions are necessary  in the remaining part of this section  and will be used freely   in the next sections.

\begin{definition}~\label{fatat}  A ring  $R$ is
\begin{enumerate}[(a)]
\item {\it reduced} if it has no non-zero nilpotent elements;
\item {\it reversible} if $ab = 0$ implies $ba = 0$ for any $a,b\in R$;
\item  said to have  {\it Insertion-of-Factors-Property} (IFP) if for $a,b\in R, ab = 0$ implies that $arb = 0$ for every $r\in R$.
\end{enumerate}
\end{definition}

\begin{definition}~\label{thatisc} An  $R$-module $M$ is
\begin{enumerate}[(a)]
 \item {\it reduced} if  whenever $a\in R$ and $m\in M$ satisfy $ma^2=0$, then $mRa =0$;
 \item  {\it symmetric}  if whenever $a,b\in R$ and $m\in M$ satisfy $mba = 0$, we have $mab = 0$;
\item said to possess IFP  if    whenever  $a\in R$ and  $m\in M$  satisfy  $ma=0$, then  $mra =0$ for  each  element  $r$  of  $R$.
\end{enumerate}
\end{definition}

\noindent The notions in the Definitions~\ref{fatat} and \ref{thatisc} have been widely studied in  \cite{buhphang2013rigid,  groenewald20132, kyomuhangi2020locally, lee2004reduced, marks2003taxonomy}.    A module   $M_R$  is said to be {\it rigid} \cite{buhphang2013rigid} if  given $a\in R$ and $m\in M$, the condition $ma^2= 0$ implies $ma = 0$.  This is equivalent to  $l_M(a^n)=l_M(a)$ for every  $a\in R$ and $n\in \Bbb Z^+$.  For commutative rings $R$, it was shown in \cite{kyomuhangi2020locally} that  $M$ is {\it reduced}  if and only if $l_M(a^n)=l_M(a)$ for every  $a\in R$ and $n\in \Bbb Z^+$.     As a dual notion to reduced modules in \cite{kyomuhangi2020locally}, we have co-reduced modules.

\begin{definition}~\label{tgat}
Let $R$ be a commutative ring.     An $R$-module  $M$ is said to be  {\it co-reduced}  if $Ma=Ma^n$ for every $a\in R$ and $n\in \Bbb Z^+$.
\end{definition}

\noindent For noncommutative rings, we give characterizations of reduced modules and reduced rings.

\begin{lemma}~\label{233} Let $R$ be a ring and $M$ be a nontrivial $R$-module. The following statements are  equivalent:
\begin{enumerate}
\item[\emph{(1)}] $M$ is reduced;
\item[\emph{(2)}] $M$ is symmetric and $l_M(a^n)=l_M(a)$  for every $a\in R$ and $n\in \Bbb Z^+$;
\item[\emph{(3)}] $M$  has IFP and $l_M(a^n)=l_M(a)$   for every $a\in R$ and $n\in \Bbb Z^+$.
\end{enumerate}
\end{lemma}
\begin{proof}~
\begin{itemize}
\item[(1)$\Leftrightarrow$(2)]  Assume that (1) holds.   By \cite[Theorem 2.2]{groenewald20132}, reduced modules are symmetric.  To prove that $l_M(a^n)=l_M(a)$, let $x\in l_M(a^n)$.  Then $xa^n=0$.  As $M$ is reduced, $xRa=0$ and so $xa=0$.  This gives $l_M(a^n)\subseteq l_M(a)$.  Since the reverse inclusion is trivial, we obtain $l_M(a^n)=l_M(a)$. The proof of $(2)\Rightarrow(1)$  holds  after applying \cite[Corollary 2.2]{groenewald20132}.
\item[(3)$\Leftrightarrow$(1)] This follows from \cite[Proposition 2.8]{buhphang2013rigid}.
\end{itemize}
\end{proof}

\begin{corollary}~\label{gag} The following statements are  equivalent for a ring $R$:
\begin{enumerate}
\item[\emph{(1)}] $R$ is reduced,
\item[\emph{(2)}] $l_R(a^n)=l_R(a)$   for every $a\in R$ and $n\in \Bbb Z^+$,
\item[\emph{(3)}] $r_R(a^n)=r_R(a)$  for every $a\in R$ and $n\in \Bbb Z^+$.
\end{enumerate}
\end{corollary}
\begin{proof}
(1)$\Leftrightarrow$(2) Reduced rings are reversible and hence have IFP.  By  Lemma~\ref{233}, (2) follows from (1).  Conversely, let $a^n=0$.   Then   $1_R\in l_R(a^n)=l_R(a)$,  so  $a=1_R\cdot a=0$.  This proves that $R$ is reduced.   The proof of (1)$\Leftrightarrow$(3) is similar.
\end{proof}

\noindent If $M_R$ is a reduced module over a commutative ring $R$, then $Ma\cong Ma^n$ for each $a\in R$ and $n\in \Bbb Z^+$.
 To see this, assume  that $M$ is reduced.    By Lemma~\ref{233}, $l_M(a^n)=l_M(a)$, and so  $Ma\cong M/l_M(a)=M/l_M(a^n)\cong Ma^n$.
For a not necessarily commutative ring $R$, the map $\varphi:M\to M$ given by $m\mapsto ma$ for $a\in R$ need not be an endomorphism.   We show in Proposition~\ref{29} that when $M$ is reduced and $e$ is an idempotent element of $R$, then $m\mapsto me$ is an idempotent endomorphism of $M$.
\begin{proposition}~\label{29}
Let $R$ be a ring,  $M$ be a reduced $R$-module, $m\in M$ and $e^2=e\in R$.  Then every map $\varphi_e$ defined by $\varphi_e(m)=me$ is an idempotent element of  $S$.
\end{proposition}
\begin{proof}
Let $e$ be an idempotent  element in $R$ and $m\in M$.   Since $M$ is reduced, it has IFP and so $me(1_R-e)=0$ implies that for any $r\in R, mer(1_R- e) = 0$, that is, $mer = mere$.  On the other hand, $m(1_R-e)e = 0$ implies that  $mre = mere$.   Hence $mer = mre$.    Now for all  $r\in R, \varphi_e(mr)=(mr)e=  (me)r=\varphi_e(m)r$.  Closure under addition always holds.
\end{proof}

\section{Weakly-endoregular modules}

\noindent Lee, Rizvi \& Roman  \cite{lee2013modules} and Ware \cite[Corollary 3.2]{ware1971endomorphism} call $M$  endoregular  if $S:=\text{End}_R(M)$ is a regular ring.
To study the various  regularity properties of the rings of endomorphisms, Anderson
and Juett   \cite{anderson2021endoregular} defined  weakly-endoregular modules.  A module $M$ over a commutative ring $R$ is  weakly-endoregular  if and only if for each $a\in R, M=Ma\bigoplus l_M(a)$.
 We give a  characterization of weakly-endoregular modules in terms of  weakly-morphic and (co-)reduced modules.     For other equivalent statements of Theorem~\ref{35} see \cite[Theorem 1.1]{anderson2021endoregular}.

\begin{theorem}~\label{35}
Let $R$ be a commutative ring and $M$ be a nontrivial $R$-module.  The following statements are  equivalent:
\begin{enumerate}
\item[\emph{(1)}] $M$ is  weakly-morphic and reduced,
\item[\emph{(2)}] $M$ is  weakly-morphic and co-reduced,
\item[\emph{(3)}] $M$ is co-reduced and reduced,
\item[\emph{(4)}] $M$ is weakly-endoregular.
\end{enumerate}
\end{theorem}
\begin{proof}~
\begin{itemize}
  \item [(1)$\Rightarrow$(2)] Assume  that (1) holds.  We need to show that  $Ma=Ma^n$ for every  $a\in R$ and $n\in \Bbb Z^+$.
 Since $M$ is weakly-morphic, $M/Ma\cong l_M(a)$ and $ M/Ma^n\cong l_M(a^n)$.
 By Lemma~\ref{233},  $l_M(a)=l_M(a^n)$, and so   $M/Ma\cong M/Ma^n$.   Therefore
 there exists an isomorphism $\varphi$ such that $\varphi(M/Ma^n)=M/Ma$.  This, $\varphi(Ma/Ma^n) = \varphi(M/Ma^n)a = Ma/Ma = 0$.   So $Ma^n = Ma$, as desired.

 \item [(2)$\Rightarrow$(1)]Assume  that (2) holds.  Then   $Ma=Ma^n$ for every  $a\in R$ and $n\in \Bbb Z^+$.   Since $M$ is weakly-morphic, $l_M(a)\cong M/Ma=M/Ma^n\cong l_M(a^n)$.   This gives $l_M(a)\cong l_M(a^n)$.   In view of \cite[Proposition 2.3]{kyomuhangi2020locally}, it remains to prove that this is equality.   By \cite[Lemma 1]{kimuli2021},  there exists some $\varphi\in S$  such that $l_M(a^n)=\varphi(M)$ and $\ker(\varphi)=Ma^n$.   The two equalities imply that $0=\varphi(M)a^n=\varphi(Ma^n)$.   By hypothesis, $Ma=Ma^n$.   So $\varphi(Ma^n) =\varphi(Ma)=\varphi(M)a$, from which  we deduce that    $l_M(a^n)=\varphi(M)\subseteq l_M(a).$  Hence $l_M(a^n)=l_M(a)$.

  \item [(2)$\Rightarrow$(3)]  Follows from the proof of $(2)\Rightarrow(1)$.

  \item [(3)$\Rightarrow$(4)]  Let $a\in R$ and assume (3).  Then $Ma=Ma^2$ by Definition~\ref{tgat} and    $l_M(a)=l_M(a^2)$ by Lemma~\ref{233}.   It follows from \cite[Theorem 1.1]{anderson2021endoregular} that  $M$ is weakly-endoregular.

  \item [(4)$\Rightarrow$(1)]  Since $M = Ma\bigoplus l_M(a),   M/Ma\cong l_M(a)$ for every $a\in R$.   Thus $M$ is weakly-morphic.    Next, let $x\in M=Ma\bigoplus l_M(a)$ such that $xa^2=0$.  Then $(xa)a=0$ and so $xa\in l_M(a)$.   But also  $xa\in Ma$ and  $Ma\cap l_M(a)=0$.      Therefore, $xa = 0$ and  $M$ is a reduced module.
\end{itemize}
\end{proof}

\begin{corollary}~\label{342}
Let $R$ be a commutative ring and $M$ be a nontrivial finitely generated $R$-module.  Then the following statements are  equivalent:
\begin{enumerate}
\item[\emph{(1)}]$M$ is  weakly-morphic and reduced,
\item[\emph{(2)}]$M$ is co-reduced,
\item[\emph{(3)}]$R/\text{Ann}_R(M)$ is a regular ring,
\item[\emph{(4)}] $M$ is weakly-endoregular.
\end{enumerate}
\end{corollary}

\begin{proof}~
\begin{itemize}
\item [(1)$\Rightarrow$(2)] Follows from Theorem~\ref{35}.
\item [(2)$\Rightarrow$(3)] Suppose that $Ma=Ma^n$ for each $a\in R$ and $n\in \Bbb Z^+$.    Then $M(a)= M(a)(a^n)$ with $M(a)$ a finitely generated module.   Using \cite[Corollary 2.5]{Atiyah1969}, we have $M(a)(1 + (a^n))=0$, which  implies that $M(a)(1+a^nr)=0$   for all  $r\in R$.  It then follows that  $(a+a^{n+1}r)\in\text{Ann}_R(M)$ and hence $a+\text{Ann}_R(M)=a^{n+1}(-r)+\text{Ann}_R(M)\in (a^{n+1})+\text{Ann}_R(M)$.   This gives $(a)+\text{Ann}_R(M)\subseteq (a^{n+1})+\text{Ann}_R(M)$ and, consequently, $(a)+\text{Ann}_R(M)=(a^{n+1})+\text{Ann}_R(M)$.   Since  $a=ra^2 + s$ for some $r\in R$ and $s\in \text{Ann}_R(M), a-ra^2\in \text{Ann}_R(M)$ and $\overline{a}=\overline{ra^2}$ for some $\overline{r}\in \overline{R}:=R/\text{Ann}_R(M)$, we have $R/\text{Ann}_R(M)$ regular.

\item [(3)$\Rightarrow$(4)] Assume that (3) holds.  Using this assumption and the First Isomorphism Theorem for the $R$-endomorphism $\varphi:R\to S:=\text{End}_R(M)$ defined by $\varphi(a)=\varphi_a$ for all $a\in R$, we obtain $\{\varphi_a:a\in R\}$ is regular.   By \cite[Proposition 7]{kimuli2021},  $M$ is a weakly-endoregular module.
\item [(4)$\Rightarrow$(1)] Follows from Theorem~\ref{35}.
\end{itemize}
\end{proof}

\begin{corollary}~\label{goad}
The following statements are  equivalent for a commutative ring $R$:
\begin{enumerate}
\item[\emph{(1)}] $R$ is  morphic and reduced,
\item[\emph{(2)}] $R$ is co-reduced,
\item[\emph{(3)}] $R$ is regular,
\item[\emph{(4)}] $R= (a)\bigoplus r_R(a)=(a)\bigoplus l_R(a)$ for each $a\in R$.
\end{enumerate}
\end{corollary}

\begin{proof}
(2)$\Leftrightarrow$(3)  A commutative ring $R$ is regular if and only if for each $a\in R,  aR=a^2R$.    Thus $R$ is co-reduced if and only if it is regular.   The  equivalences (1)$\Leftrightarrow$(3)$\Leftrightarrow$(4) follow from Theorem~\ref{35}.
\end{proof}

\begin{corollary}~\label{}
Every module over a commutative regular ring is weakly-endoregular.
\end{corollary}

\begin{proof}
This follows from  \cite[Proposition 13]{kimuli2021}.
\end{proof}

\section{Abelian endoregular modules}
\noindent The focus of  this section is the characterization  of Abelian endoregular modules in terms of reduced  and morphic modules.
A ring $R$ is said to be {\it Abelian} if all its idempotents are central.  If $R$ is reduced, then every idempotent is central.    A strongly regular ring is reduced, regular and Abelian.  More generally, an $R$-module $M$ is said to be  {\it  Abelian} if $S$ is an Abelian ring.
 $M_R$ is  an  {\it Abelian endoregular module} if  $S$ is a regular and Abelian ring.

\begin{remark}~\label{30}
$M_R$ is  an  Abelian endoregular module  if and only if  $M=\varphi(M)\bigoplus\ker(\varphi)$ for every $\varphi\in S$.  Abelian endoregular modules are morphic modules.    Note that an endoregular module need not be morphic. The  $\Bbb Z$-module $\bigoplus_{i=1}^\infty\Bbb Q_i$, where $\Bbb Q_i=\Bbb Q$,   is endoregular but not morphic.
\end{remark}

\noindent Recall that $M$ {\it cogenerates}  $M/\varphi(M),\varphi\in S$ if $M/\varphi(M)$  can be embedded in $M^{(I)}$, where $I$ is an index set. That is,
 $0\neq x\in M/\varphi(M),\varphi\in S$, implies that    $\gamma(x)\neq 0$  for some  $\gamma\in\text{Hom}_R(M/\varphi(M),M)$ \cite[pg. 230]{nicholson1993pp}.

\begin{lemma}~\label{27}
If $R$ is a ring and $M$ a nontrivial morphic $R$-module, then  for every $\varphi\in S,$ $$\varphi(M)=r_M(l_S(\varphi)).$$
Moreover, the following statements are equivalent:
\begin{enumerate}
\item[\emph{(1)}] $\varphi(M)=r_M(l_S(\varphi))$ for every $\varphi\in S$;
\item[\emph{(2)}] For each $\varphi\in S$ and $m\in M$, if $l_S(\varphi(M))\subseteq l_S(m)$, then $m\in \varphi(M)$;
\item[\emph{(3)}]$M$ cogenerates $M/\varphi(M)$ for each $\varphi\in S$.
\end{enumerate}
\end{lemma}
\begin{proof}
For every $\varphi\in S$, there exists $\psi\in S$ such that $\varphi(M)=\ker(\psi)=r_M(\psi)$ and $\psi(M)=\ker(\varphi)=r_M(\varphi)$.  It follows from the equality   $\varphi(M)=r_M(\psi)$
 that $\psi\varphi(M)=0$ which gives $\psi\varphi=0$ and hence
 $S\psi\subseteq l_S(\varphi(M))$.
 Thus  $r_M(l_S(\varphi))\subseteq r_M(\psi)=\varphi(M)$.   The reverse inclusion is obvious, hence $r_M(l_S(\varphi))=\varphi(M)$.
\begin{itemize}
\item [(1)$\Rightarrow$(2)]  Let $m\in M$ and $\varphi\in S$ such that $l_S(\varphi(M))\subseteq l_S(m)$.    Then $m\in r_M(l_S(m))\subseteq r_M(l_S(\varphi(M)))=\varphi(M)$ by (1).  Hence $m\in \varphi(M)$.
\item [(2)$\Rightarrow$(1)]Clearly  $\varphi(M)\subseteq r_M(l_S(\varphi(M)))$ for every $\varphi\in S$.  Let $m\in r_M(l_S(\varphi(M)))$.  Then  we have $l_S(r_M(l_S(\varphi(M))))\subseteq l_S(m)$.  By \cite[Proposition 24.3]{anderson1992rings}, $l_S(\varphi(M))\subseteq l_S(m)$.  By (2), $m\in \varphi(M)$.
\item [(1)$\Leftrightarrow$(3)] Assume that (1) holds and let $\varphi\in S$.  In view of \cite[pg. 109 and Lemma 24.4]{anderson1992rings},   $\text{Rej}_{M/\varphi(M)}(M):=\bigcap\{\ker(\gamma):\gamma\in \text{Hom}_R(M/\varphi(M),M)\}=r_M(l_S(\varphi(M)))/\varphi(M)=0$.   This gives $\text{Rej}_{M/\varphi(M)}(M)=0$ for each $\varphi\in S$.   Applying
       \cite[Corollary 8.13]{anderson1992rings}, $M$ cogenerates $M/\varphi(M)$ for each $\varphi\in S$.   Conversely, suppose $M$ cogenerates $M/\varphi(M)$ for each $\varphi\in S$.  Then $\text{Rej}_{M/\varphi(M)}(M)=0$ for each $\varphi\in S$ by \cite[Corollary 8.13]{anderson1992rings}.   Applying \cite[Lemma 24.4]{anderson1992rings}  gives $r_M(l_S(\varphi(M)))/\varphi(M)=0$.     Hence $r_M(l_S(\varphi(M)))=\varphi(M)$ follows.
\end{itemize}
\end{proof}

\begin{remark}~\label{waki}
\begin{enumerate}[(a)]
\item In view of Lemma~\ref{27}, the hypothesis ``$\varphi(M)=r_M(l_S(\varphi(M)))$'' in statement (b) of \cite[Proposition 4.2]{lee2013modules} is superfluous.
\item  Recall that in \cite{nicholson1993pp}, a module $_{S}M$ is  $P$-{\it injective} if  $\varphi(M)=r_M(l_S(\varphi))$ for every $\varphi\in S$.  A ring $R$ is called {\it left $P$-injective} if it is a $P$-injective right $R$-module (equivalently, $r_Rl_R(a)=aR$ for every $a\in R$).   Thus, if $M$ is a morphic $R$-module, then  ${}_SM$ is a $P$-injective module by Lemma~\ref{27}.
\end{enumerate}
\end{remark}

\begin{lemma}~\label{topp}
If $M$ is a morphic $R$-module and $S$ is a reduced ring,  then  for every $\varphi\in S$,
 $$r_M(\varphi^2)=r_M(\varphi).$$
\end{lemma}
\begin{proof}
We only prove  $r_M(\varphi^2)\subseteq r_M(\varphi)$ since the reverse inclusion is obvious.        Since $M$ is morphic, there exists $\gamma\in S$ such that $\gamma(M)=r_M(\varphi^2)$.  This implies $\varphi^2\gamma=0$.   Further, $S$ being reduced implies that $\varphi\gamma=0$.  So, $\gamma(M)\subseteq r_M(\varphi)$ and we get $r_M(\varphi^2)=\gamma(M)\subseteq r_M(\varphi)$.
\end{proof}

\noindent Now we give a characterization of Abelian endoregular modules in terms of morphic modules and reduced rings of endomorphisms.
\begin{theorem}~\label{31}
Let $R$ be a ring and $M$ be  a nontrivial $R$-module.  The following statements are equivalent:
\begin{enumerate}
\item[\emph{(1)}] $M_R$ is a morphic module  and $S$ is a reduced ring,
\item[\emph{(2)}] $M_R$ is an Abelian endoregular module.
\end{enumerate}
\end{theorem}

\begin{proof}~
\begin{itemize}
\item [(1)$\Rightarrow$(2)]  Since $S$ is reduced,     $l_S(\varphi)=l_S(\varphi^2)$ for any   $\varphi\in S$ by Corollary~\ref{gag}.    Applying Lemma~\ref{27}, $\varphi(M)=r_M(l_S(\varphi(M)))=r_M(l_S(\varphi^2(M)))=\varphi^2(M)$.  This gives $\varphi(M)=\varphi^2(M)$.
 For any $m\in M, \varphi(m)\in \varphi(M)=\varphi^2(M)$, and   so  $\varphi(m) = \varphi(n)$ for some $n\in \varphi(M)$.    Therefore, $x: = m-n\in r_M(\varphi)$ and $m= x + n\in r_M(\varphi) + \varphi(M)$.   We obtain  $M = r_M(\varphi) + \varphi(M)$. To prove that  this is a direct sum, let $y\in r_M(\varphi)\cap\varphi(M)$.   Then $y=\varphi(m)$ for some $m\in M$  with $\varphi(y)=\varphi^2(m)=0$.
Consequently we have $m\in r_M(\varphi^2)=r_M(\varphi)$ by Lemma~\ref{topp}, from which we have $y=\varphi(m)=0$, thus we obtain  $0=r_M(\varphi)\cap\varphi(M)$.   This proves $M = r_M(\varphi)\bigoplus\varphi(M)$, and  $M$ is an Abelian endoregular module.
\item [(2)$\Rightarrow$(1)]  $M$ is morphic by  Remark~\ref{30}.  In addition,   since  $S$ is a strongly regular ring, it is reduced.
\end{itemize}
\end{proof}

\begin{corollary}~\label{thats}
Let $R$ be a ring and $M$  a nontrivial $R$-module. The following statements are equivalent:
\begin{enumerate}
\item[\emph{(1)}] $M_R$ is  morphic  and ${}_SM$ is  reduced,
\item[\emph{(2)}] $M_R$ is an Abelian endoregular module.
\end{enumerate}
\end{corollary}
\begin{proof}~
\begin{itemize}
\item [(1)$\Rightarrow$(2)]  Let $\varphi\in S$  such that $\varphi^2=0$.    Then $\varphi S(m)=0$ for all $m\in M$ by (1).     It follows that   $\varphi 1_M(m) = \varphi(m) = 0$ for all  $m\in M$, so $\varphi = 0$.  This shows that $S$ is a reduced ring and  (2) follows  by Theorem~\ref{31}.

\item [(2)$\Rightarrow$(1)] $M_R$ is clearly morphic by Remark~\ref{30}.  Let   $\varphi\in S$ and $m\in M$ such $\varphi^2(m)=0$.   Then $\varphi (\varphi(m))=0$.   In view of  Remark~\ref{30},  $\varphi(m)\in r_M(\varphi)\cap\varphi(M)=0$, so $\varphi(m)=0$.
 Since $S$ is strongly regular, there exists some $\psi\in S$ such that $\varphi=\varphi\psi\varphi$ with $\varphi\psi=\psi\varphi$   a central idempotent element  of $S$. Thus $\varphi S(m)=\varphi\psi\varphi S(m)=\varphi S\psi\varphi(m)=0$.   This proves that ${}_SM$ is a  reduced module.
 \end{itemize}
\end{proof}

\noindent By considering the case $M=R$ and $\text{End}_R(R)\cong R$, we have:

\begin{corollary}~\label{336}
The following statements are equivalent for a  ring $R$:
\begin{enumerate}
\item[\emph{(1)}] $R$ is right morphic and reduced,
\item[\emph{(2)}] $R$ is left $P$-injective and reduced,
\item[\emph{(3)}]  $R= aR\bigoplus r_R(a)$ for each $a\in R$,
\item[\emph{(4)}] $R$ is  strongly regular.
\end{enumerate}
\end{corollary}
\begin{proof}~
\begin{itemize}
\item [(1)$\Rightarrow$(2)]  This follows by Lemma~\ref{27} and Remark~\ref{waki} (b).
\item [(1)$\Leftrightarrow$(3)] and  (3)$\Leftrightarrow$(4)  These are a consequence  of  Theorem~\ref{31} and Remark~\ref{30}.
\item [(2)$\Rightarrow$(4)]   Since  $R$ is reduced, for each $a\in R, l_R(a) = l_R(a^2)$ by Corollary~\ref{gag}.  It follows that $aR=r_R(l_R(a))=
r_R(l_R(a^2))= a^2R$ because $R$ is left $P$-injective.  Thus $a=a^2y$ for some $y\in R$, and this proves $R$ is strongly regular.
\end{itemize}
\end{proof}

\noindent A module $M_R$ is  {\it duo} provided every   sub-module of $M$ is fully invariant, that is, for any sub-module $N$ of $M,\varphi(N)\subseteq N$ for every $\varphi\in S$.   A ring  $R$ is {\it right duo} if and only if every right ideal of $R$ is a two-sided ideal,  equivalently if $Ra$ is contained in $aR$ for every element $a$ in $R$ \cite{ozcan2006duo}.

\begin{lemma}{\rm \cite[Lemma 1.1]{ozcan2006duo}}~\label{6} Let $R$ be any ring.  Then a right $R$-module $M$ is a duo module if and only if for each endomorphism $\varphi$ of $M$ and each element $m$ of $M$ there exists $a$ in $R$ such that $\varphi(m) = ma$.
\end{lemma}

\begin{lemma}~\label{26} Let $R$ be a commutative ring.  For  a nontrivial duo $R$-module $M$,  consider the following statements:
\begin{enumerate}
\item[\emph{(1)}]$M$ is reduced as a right $R$-module,
\item[\emph{(2)}] $M$ is  reduced as a left $S$-module,
\item[\emph{(3)}] $S$ is a reduced ring.
\end{enumerate}
Then $(1)\Leftrightarrow(2)\Rightarrow(3)$.
\end{lemma}
\begin{proof}~
\begin{itemize}
\item [(1)$\Rightarrow$(2)]  Let $\varphi\in S$ and $m\in M$ such that $\varphi^2(m)=0$.  Then    $ma^2=0$ for some $a\in R$ because $M$ is duo.  By (1),    $mra=0$ for all $r\in R$.    Since  every element  in $S$ is defined by right multiplication of each element of $M$ by some element of $R, \varphi(\psi(m))=mra=0$ for every $\psi\in S$ for some $r\in R$. Thus $\varphi S(m)=0$ and ${}_SM$ is a reduced module.

 \item [(2)$\Rightarrow$(1)] Let $m\in M$ and $a\in R$ such that $ma^2=0$.   Then the endomorphism  $\varphi:M\to M, x\mapsto xa$ gives $\varphi^2(m)=0$.
  Since ${}_SM$ is reduced, we have $\varphi S(m)= 0$.   Note that for every $r\in R$, right multiplication by $r$ defines an endomorphism $\psi_r:M\to M,m\mapsto mr$.   This gives $mra=\varphi\psi_r(m)=0$.  Since
   $mRa\subseteq\varphi S(m)=0,   M_R$ is a reduced module.

\item [(2)$\Rightarrow$(3)]  Let $\varphi\in S$  such that $\varphi^2=0$.   Then for each $m\in M$, Lemma~\ref{6} gives $0=ma^2=\varphi^2(m)$ for some $a\in R$; and  so $\varphi S(m)=0$ by (2).  It follows that   $\varphi 1_M(m) = \varphi(m) = 0$.   Since $m$ was chosen arbitrarily,   $\varphi = 0$.

  \end{itemize}
\end{proof}

\noindent Note that even when  a duo module has a reduced  ring of endomorphisms, the module itself may not be reduced.
 \begin{Ex}~\label{tred}
Let $R:=\Bbb Z$ be a ring.   For any prime $p$, the Pr$\ddot{\text{u}}$fer $p$-group $M:=\Bbb Z(p^\infty)$ is an Artinian uniserial $R$-module and hence  a duo module by {\rm\cite[pg. 536]{ozcan2006duo}}.   Then it is well-known that $S:=\text{End}_{R}(M)$ is the ring of $p$-adic integers {\rm\cite[Exercises 3 (17), pg. 54]{anderson1992rings}}.   Since the ring of $p$-adic integers is a commutative domain, it is a reduced ring.   However, $M$ is neither reduced as an $R$-module nor as an  $S$-module.
\end{Ex}

\noindent An $R$-module $M$ is said to be a {\it multiplication module} provided for each submodule $N$ of $M$ there exists an ideal $A$ of $R$ such that $N =MA$.  Finitely generated multiplication modules are duo.

\begin{lemma}{\rm\cite[Proposition 19]{kimuli2021}}~\label{300}
Let $M$ be a finitely generated multiplication module over a commutative ring $R$.  Then  $M$  is weakly-morphic if and only if it is morphic.
\end{lemma}

\begin{corollary}~\label{tgt}
Every cyclic module  over  a commutative ring $R$  is weakly-morphic if and only if it is morphic.
\end{corollary}
\begin{proof}
Since every cyclic $R$-module is a multiplication module that is finitely
generated, it is weakly-morphic if and only if it is morphic by Lemma~\ref{300}.
\end{proof}

\begin{proposition}~\label{444}
Let $R$ be a commutative ring and $M$ be a nontrivial finitely generated multiplication $R$-module.   Then $M$ is weakly-endoregular if and only if it is Abelian endoregular.
\end{proposition}
\begin{proof}
Assume that $M_R$ is a weakly-endoregular module.
  By Theorem~\ref{35}, $M_R$ is  weakly-morphic and reduced.
As it is a finitely generated multiplication module, Lemma~\ref{300} implies  $M_R$ is    morphic.
    Being duo, $_{S}M$ is a reduced module by Lemma~\ref{26}.
     Applying Corollary~\ref{thats} proves that $M$ is Abelian endoregular.   The converse clearly holds since every Abelian endoregular module over a commutative ring is weakly-endoregular.
\end{proof}

\begin{corollary}~\label{unu}
Every cyclic module over  a commutative ring $R$  is weakly-endoregular if and only if it is Abelian endoregular.
\end{corollary}

\begin{proof}
Since cyclic modules are finitely generated multiplication modules, the proof of the corollary is immediate from Proposition~\ref{444}.
\end{proof}

\noindent  An $R$-module $M$ is {\it strongly duo} \cite{khabazian2010strongly} if the {\it trace} of $M$ in $N$ is $N$, that is, $\text{Tr}_N(M):=\sum\{\text{Im}(\lambda):\lambda\in\text{Hom}_R(M,N)\}= N$ for all $N\subseteq M_R$.  Clearly, every strongly duo module $M$ is a duo module.   In  \cite[Theorem 5.5]{khabazian2010strongly}, the  ring of endomorphisms of  a module $M$ that is  strongly duo and reduced  was shown to be a strongly regular ring.    For commutative rings, we have an improved result in Corollary~\ref{32}.

\begin{corollary}~\label{32}
Let $R$ be a commutative ring and $M$ be a  nontrivial duo $R$-module.  The following statements are equivalent:
\begin{enumerate}
\item[\emph{(1)}] $M_R$ is a morphic and reduced module,
\item[\emph{(2)}]$S$ is a strongly regular ring.
\end{enumerate}
\end{corollary}
\begin{proof}~
\begin{itemize}
\item [(1)$\Rightarrow$(2)] Assume (1) holds.  Then  $S$ is a reduced ring by  Lemma~\ref{26} and is, therefore, strongly regular by Theorem~\ref{31}.
\item [(2)$\Rightarrow$(1)]  $M_R$ is  morphic by Remark~\ref{30}  and    reduced by    Lemma~\ref{26}.
\end{itemize}
\end{proof}

\section{F-regular modules}

\noindent  Recall that a ring $R$ is  regular if and only if every right (left) cyclic ideal of $R$ is a direct summand of $R_R$.   To generalise this  characterization to modules, Ramamurthi and  Rangaswamy  in \cite[pg. 246]{ramamurthi1973finitely} defined  strongly regular modules.   A module $M$ is called {\it strongly regular} (in the sense of \cite{ramamurthi1973finitely}) if every finitely generated submodule is a direct summand, or equivalently every  cyclic submodule is a direct summand.       Following  Naoum \cite{naoum1995regular},   we call the strongly regular modules {\it strongly F-regular} (even without commutativity of $R$).
In \cite{nicholson2005survey}, a relationship between morphic finitely generated strongly F-regular modules and their  rings of endomorphisms  was established.

\begin{proposition}{\rm\cite[Corollary 2.7]{nicholson2005survey}}~\label{7}
A finitely generated strongly F-regular module $M$ is morphic if and only if $S$ is morphic and regular.
\end{proposition}

\noindent An $R$-module $M$ is said to be {\it $k$-local-retractable}  (for kernel-local-retractability) (or equivalently,  {\it P-flat over $S:=\text{End}_R(M)$}) if  for any $\varphi\in S$ and any nonzero element $x\in r_M(\varphi)$, there exists a homomorphism $\psi_x:M\to r_M(\varphi)$ such that $x\in\psi_x(M)\subseteq r_M(\varphi)$  (\cite[pg. 4069]{lee2013modules} and \cite{nicholson1993pp}).
The module $M_R$ is called a {\it self-generator} in \cite[pg. 228]{nicholson1993pp} if it generates each of
its images, that is, $mR= \text{Hom}_R(M,mR)(M)$ for all $m\in M$.  In this case, for each $m\in M,m=\sum \alpha_i(x_i)$ with $x_i\in M$ and $\alpha_i\in \text{Hom}_R(M,mR)$.
\begin{proposition}~\label{wewer}
Every nontrivial strongly F-regular module $M_R$ is a k-local-retractable module.
\end{proposition}
\begin{proof}
Since strongly F-regular modules are self-generator modules by \cite[pg. 228]{nicholson1993pp}, $M_R$ is P-flat over $S$ by \cite[Lemma 1]{nicholson1993pp}, which is equivalent to being
$k$-local-retractable  by  \cite[pg. 4069]{lee2013modules}.
\end{proof}
\begin{lemma}~\label{ewaffe}
If $M$ is a $k$-local-retractable $R$-module and  $S$ is a reduced ring, then for every $\varphi\in S$,
 $$r_M(\varphi^2)=r_M(\varphi).$$
\end{lemma}
\begin{proof}
Let $x\in r_M(\varphi^2)$.   Due to $k$-local-retractability of $M$, there exists $0\neq \psi_x\in S$ such that $x\in \psi_x(M)\subseteq r_M(\varphi^2)$.   Hence  $\varphi^2\psi_x= 0$.   Since $S$ being reduced implies $\varphi\psi_x= 0, x\in\psi_x(M)\subseteq r_M(\varphi)$.  This shows that $r_M(\varphi^2)\subseteq r_M(\varphi)$.  The reverse inclusion is well-known.
\end{proof}

\begin{lemma}~\label{ewafee}
If  $M$ is a nontrivial duo and  strongly F-regular $R$-module, then $S$ is a reduced ring.
\end{lemma}
\begin{proof}
Let $\varphi\in S$ such that $\varphi^2=0$.  If $\varphi\neq 0$, then there exists some $0\neq m\in M$ such $\varphi(m)\neq0$.  By the strongly F-regular hypothesis, $M=\varphi(m)R\bigoplus X$ for some submodule $X$ of $M$.  Since $M$ is duo, $\varphi(M)=\varphi(\varphi(m)R)\bigoplus\varphi(X)=\varphi(X)\subseteq X$, so $\varphi(M)\subseteq X$.  This implies that $\varphi(m)\in \varphi(m)R\cap X=0$, a contradiction.  Thus $\varphi=0$ and $S$ is a reduced ring.
\end{proof}

\begin{lemma}~\label{ewafe}
Let $M$ be a nontrivial duo and  strongly F-regular $R$-module.
  If $K\cong K'$ where $K$ and $K'$ are  submodules of $M$, then   $K=K'$.
\end{lemma}
\begin{proof}
First, we  prove that for every submodule $N$ of $M,  \varphi(N) \subseteq N$ for all homomorphisms  $\varphi:N\to M$.
Let  $n\in N$ and consider $\varphi:nR\to M$.      By the strongly F-regular hypothesis, $M=nR\bigoplus X$ for some submodule $X$.     Define $\beta:M\to M$ by $\beta(s+x)=\varphi(s)$ for every  $s\in nR$ and $x\in X$.   Then $\beta$  is a well-defined endomorphism of $M$ which   extends $\varphi$ to an endomorphism
of    $M$.    It follows that for any   $n\in N$ there exists $\beta\in S$ such that $\varphi(n)\in\varphi(nR)=\beta(nR)\subseteq N$ because $M$ is duo.   Hence $\varphi(N)\subseteq N$.
Therefore, if $\sigma:K\to K'$ is the given isomorphism, then $K'=\sigma(K)\subseteq K$  and $K=\sigma^{-1}(K')\subseteq K'$.  This proves that $K=K'$.
\end{proof}

\noindent  The following equivalent conditions were established in \cite[Proposition 4.13 and Lemma 4.2]{bamunoba2020morphic} for near-rings, so they must hold for rings:  reduced and right morphic $\Leftrightarrow$  regular and right duo $\Leftrightarrow$  reduced and regular $\Leftrightarrow$   strongly regular.
 In the next theorem we write down these ideas in the module-theoretic context.

\begin{theorem}~\label{3} Let $R$ be a ring and  $M$ be a nontrivial strongly F-regular module.  Then the following statements are equivalent:
\begin{enumerate}
\item[\emph{(1)}] $M_R$ is a morphic module and $S$ is a reduced ring,
\item[\emph{(2)}]  $M_R$ is a duo module,
\item[\emph{(3)}] $M_R$ is an Abelian endoregular module.
\end{enumerate}
\end{theorem}
\begin{proof}~
\begin{itemize}
\item[(1)$\Rightarrow$(2)]Assume that (1) holds.   Let $N$ be a submodule of $M$   and $\varphi\in S$.  By the  strongly F-regular hypothesis, for every $n\in N,nR = e(M)$ for some idempotent $e\in S$.    Since $S$ is reduced, $e$ is  central in $S$.  Hence, $\varphi(n)\in \varphi(nR)= \varphi (e(M))= e(\varphi(M))\subseteq e(M)=nR \subseteq N$.   This proves that $\varphi(N) \subseteq N$ for all $\varphi\in S$,  so $M_R$ is duo.

\item[(2)$\Rightarrow$(1)] Assume that (2) holds.  Then $S$ is  a reduced ring by Lemma~\ref{ewafee}.   To prove $M$ is morphic, in view of Theorem~\ref{31}, we will show that $M=\varphi(M)\bigoplus r_M(\varphi)$   for each $\varphi\in S$.
 Let $\varphi\in S$.   Using Lemma~\ref{ewaffe} and the First Isomorphism Theorem, $\varphi(M)\cong M/r_M(\varphi)=M/r_M(\varphi^2)\cong \varphi^2(M)$.
 This gives $\varphi(M)\cong\varphi^2(M)$.   Applying Lemma~\ref{ewafe} gives $\varphi(M)=\varphi^2(M)$.
 For any $x\in M, \varphi(x)\in \varphi(M)=\varphi^2(M)$ which implies there exists $y\in M$ such that $\varphi(x) = \varphi^2(y)$.  Then
$\varphi(x -\varphi(y))=0$.  This implies $k:=x-\varphi(y)\in r_M(\varphi)$, hence $x=\varphi(y)+k\in \varphi(M)+r_M(\varphi)$ and $M=\varphi(M)+r_M(\varphi)$.     Let $x\in r_M(\varphi)\cap\varphi(M)$.   Then $x=\varphi(m)$ for some $m\in M$  with $\varphi(x)=\varphi^2(m)=0$.
Consequently, in view of Proposition~\ref{wewer} and  Lemma~\ref{ewaffe}, we have $m\in r_M(\varphi^2)= r_M(\varphi)$, from which we have $x=\varphi(m)=0$.
 Thus $0=r_M(\varphi)\cap\varphi(M)$ and $M = \varphi(M)\bigoplus r_M(\varphi)$.
\item[(1)$\Leftrightarrow$(3)] Follows from Theorem~\ref{31}.
 \end{itemize}
\end{proof}

\begin{definition}~\label{trew}
A submodule $N$ of $M$  is {\it pure} in $M$ if the sequence $0\to N\bigotimes E\to M\bigotimes E$  is exact for each $R$-module $E$.
$N$  is  {\it relatively divisible-pure} or {\it RD-pure}  in  $M$  in case  $Na=Ma\cap N$ (equivalently, $0\to N\bigotimes R/aR\to M\bigotimes R/aR$ is exact)  for each $a\in R$.
 \end{definition}
\noindent  By \cite[Proposition 8.1]{Fieldhouse}, every pure submodule is also RD-pure.
A ring $R$ is regular if and only if every (right) ideal is pure (see \cite{Fieldhouse}).   Using this fact, Fieldhouse calls $M_R$  a {\it regular module}  if every submodule $N$ of  $M$ is pure.    Following Naoum \cite{naoum1995regular}, we call the  Fieldhouse  regular modules $F$-{\it regular}.

\begin{lemma}~\label{trey}
Let $R$ be a commutative ring and $M$ be a nontrivial  F-regular  module.  Then $M$ is  weakly-endoregular,
 weakly-morphic, reduced and co-reduced.
 \end{lemma}
\begin{proof}  Since submodules of F-regular  modules are RD-pure by \cite[Proposition 8.1]{Fieldhouse}, we  show that $M=Ma\bigoplus l_M(a)$  for each $a\in R$.      Let $a\in R$.     Then  $Ma=Ma\cap Ma=Ma^2$ so that $Ma=Ma^2.$      It follows that for any   $x\in M,   xa=na^2$ for some $n\in M$.   Since $(x-na)a=0, x-na\in l_M(a)$ and     $x=na+x-na\in Ma +l_M(a)$.
Hence $M=Ma+l_M(a)$.      By the RD-pure property, $0=l_M(a)a=Ma\cap l_M(a)$ for every $a\in R$.  Thus   $M=Ma\bigoplus l_M(a)$.
   This proves that $M$ is weakly-endoregular.   Using Theorem~\ref{35},  $M$ is weakly-morphic,    reduced and co-reduced.
\end{proof}

\begin{Ex}~\label{tgvr}
The converse of Lemma~\ref{trey} does not hold in general.
The $\Bbb Z$-module $\Bbb Q$ is weakly-endoregular, weakly-morphic and reduced  but it is  not F-regular.   In particular,  not all its submodules are  (RD-)pure since $2\Bbb Q\cap\Bbb Z\neq 2\Bbb Z$ for the submodule $\Bbb Z$.
\end{Ex}

\noindent Since submodules of strongly F-regular modules are RD-pure by \cite[pg. 240 and 246]{ramamurthi1973finitely},  it follows from \cite[Proposition 8]{kimuli2021} that if $R$ is a commutative ring, then every strongly F-regular module  is a weakly-morphic module.
An $R$-module $M$ is  {\it finitely presented} (abbreviated as f.p.)  if there exists an exact sequence of the form $R^n\to R^m\to M$ with $n,m\in \Bbb Z^+$, or equivalently if   $M\cong P/Q$, where $P$ and $Q$ are finitely generated modules, and $P$ is a projective module.  Clearly,  strongly F-regular modules are F-regular but the converse is not true in general, see \cite{anderson2019module}.  In Proposition~\ref{346},  we determine when the F-regular modules are strongly F-regular.

\begin{lemma}~{\rm\cite[Theorem 7.14]{matsumura1989commutative}}~\label{345}
If $N$ is a pure submodule of $M$ and  $M/N$  is finitely presented, then $N$ is a direct summand  of  $M$.
\end{lemma}

\begin{proposition}~\label{346}
Let $R$ be a commutative ring and $M$ be a nontrivial $R$-module.   Then $M$ is strongly F-regular  whenever  $M$ is F-regular and $M/mR$ is finitely presented  for each $m\in M$.
\end{proposition}
\begin{proof}
Suppose $M$ is F-regular and  $M/mR$ is  finitely presented  for each $m\in M$.    Then $mR$ is a pure submodule in $M$ for each $m\in M$.  By Lemma~\ref{345}, $mR$ is a direct summand of $M$  for each $m\in M$ and, thus $M$ is strongly F-regular.
\end{proof}

\noindent Let  $R$ be a commutative ring and  $N$ be a proper submodule of $M_R$.   $N$ is a {\it prime submodule} if for any $a\in R$ and $m\in M, ma\in N$ implies either $m\in N$ or $a\in (N :_R M):=\{r\in R : Mr\subseteq N\}$.  For any proper submodule $N$ of $M$, the intersection of all  prime submodules of $M$ containing $N$ is denoted by $Rad(N)$.
Theorem~\ref{230}  gives some new characterizations for F-regular  modules over commutative rings.   For other equivalent statements of Theorem~\ref{230} see  {\rm\cite[Theorem 6]{anderson2019module}}, \cite[Theorem 2.3, Corollary 2.7 and Theorem 4.1]{hassanzadeh2019regular} and  \cite[Theorem 2.1]{sharif1996rings}.

\begin{theorem}~\label{230}
Let $R$ be a commutative ring and $M$ be a nontrivial $R$-module.  The following statements are equivalent:
\begin{enumerate}
\item[\emph{(1)}] $M$ is F-regular,
\item[\emph{(2)}] Every  submodule  of $M$ is a weakly-endoregular module,
\item[\emph{(3)}] Every  submodule  of $M$ is a weakly-morphic and reduced module,
\item[\emph{(4)}] Every cyclic submodule of $M$ is a (weakly-)morphic and reduced module,
\item[\emph{(5)}] Every cyclic submodule of $M$ is a  co-reduced module,
\item[\emph{(6)}] Every cyclic submodule of $M$ is an  Abelian endoregular module,
\item[\emph{(7)}] Every cyclic submodule of $M$ is an  F-regular module.
\end{enumerate}
\end{theorem}
\begin{proof}~
\begin{itemize}
\item[(1)$\Rightarrow$(2)] Assume that  (1) holds and let $N$ be a submodule of $M$.   By Lemma~\ref{trey}, $M$ is  weakly-endoregular.  Since $N$ is pure in $M$, by \cite[Theorem 1.1 (3)]{anderson2021endoregular} $N$ is  weakly-endoregular as well.
\item[(2)$\Rightarrow$(3)]  This follows from Theorem~\ref{35}.
\item[(3)$\Rightarrow$(4)] Using Corollary~\ref{tgt}, weakly-morphic  cyclic modules are morphic.
\item[(4)$\Rightarrow$(5)] Since every cyclic submodule of $M$ is a finitely generated $R$-module, the proof follows from Corollary~\ref{342}.
\item[(5)$\Rightarrow$(1)] Suppose that (5) holds.   Let $N$ be a proper submodule of $M$.  In view of \cite[Theorem 2.3]{hassanzadeh2019regular}, we have  to prove that $Rad(N)=N$.    But to prove that $Rad(N)=N$, by \cite[Theorem 2.1]{sharif1996rings}, it is enough to show that $m(a) = m(a^2)$ for all $a\in R$ and $m\in M$.  Let  $a\in  R$ and $m\in M$.   Since $mR$ is a co-reduced module,   $mRa=mRa^2$.  Thus  $m(a) = m(a ^2 )$.
\item[(2)$\Rightarrow$(6)] Assume  (2) holds.   Since $mR,m\in M$ is a finitely generated multiplication module,  it is  weakly-endoregular  if and only if  it is Abelian endoregular by Proposition~\ref{444}.

\item[(6)$\Rightarrow$(5)] Assume that $mR$ is Abelian endoregular module  for each $m\in M$.   Then $mR=mRa\bigoplus l_{mR}(a)$ for each $m\in M$ and $a\in R$.  It follows that $mRa=mRa^2$, and $mR$ is co-reduced for each $m\in M$ by Definition~\ref{tgat}.
\item[(1)$\Rightarrow$(7)] Assume that $M$ is F-regular.   Then $mR$ is an F-regular module  for each $m\in M$ by \cite[Theorem 8.2]{Fieldhouse} and \cite[Proposition 2.6]{hassanzadeh2019regular}.
\item[(7)$\Rightarrow$(1)] Assume $mR$ is an F-regular module for each $m\in M$.   Then by Definitions~\ref{thatis} and \ref{trew}, $mRa$ is a(n) (RD-)pure submodule of $mR$ for each  $a\in R$.   It follows that $mRa=mRa\cap mRa=mRa^2$, proving that $m(a) = m(a ^2 )$.   By \cite[Theorem2.1]{sharif1996rings}, $Rad(N)=N$ for each submodule $N$ of $M$.  Hence  $M$ is F-regular by \cite[Theorem 2.3]{hassanzadeh2019regular}.
\end{itemize}
\end{proof}

\begin{remark}~
\begin{enumerate}[(a)]
\item If $M$ is an F-regular module over a commutative ring,  then  by Lemma~\ref{trey} and Theorem~\ref{230},  $\ker(\varphi)$ and $\text{Im}(\varphi)$ are  weakly-morphic and reduced modules for every $\varphi\in S$.

\item  By {\rm Theorem~\ref{230}}, the properties: ``weakly-morphic module'' and ``reduced module''  transfer from a module to each of its  submodules and conversely.
    \item   It is shown in  {\rm Theorem~\ref{230}} that if every (cyclic) submodule of $M$ is a (weakly-)morphic and reduced module, then  $M$ attains the F-regularity property.
\end{enumerate}

\end{remark}

\noindent Recall that a commutative ring $R$ is regular if and only if for each $a\in R,  aR=a^2R$.
For commutative rings $R$,  Jayaram and Tekir in \cite{jayaram2018neumann} call $M_R$  {\it regular} if for each  $m\in M,mR=Ma=Ma^2$ for some $a\in R$.  Following \cite[Definition 1]{anderson2019module}, we  call the Jayaram and Tekir regular modules  {\it JT-regular}.     Anderson,  Chun  \& Juett  in \cite{anderson2019module} defined a weak version of these modules, the   weakly JT-regular modules.   $M$ is a  {\it weakly JT-regular module} if  $Ma=Ma^2$ for each $a\in R$.

\begin{remark}~\label{231x} Let $R$ be a commutative ring and $M$ be an $R$-module.
\begin{enumerate}[(a)]
\item By \cite[Theorem 13]{anderson2019module}, $M$ is JT-regular $ \Rightarrow M$  is strongly F-regular $ \Rightarrow M$ is F-regular $ \Rightarrow M$ is weakly JT-regular.
\item By Definitions~\ref{thatis} and \ref{tgat}, the weakly JT-regular modules and the  co-reduced modules are indistinguishable.   Therefore, by (a)   F-regular (resp., strongly F-regular, JT-regular) modules  are co-reduced modules.

\item The fact that  finitely generated JT-regular modules are reduced was proved in \cite[Lemma 10]{jayaram2018neumann}.    By Lemma~\ref{trey}, if $M$ is  an F-regular (resp., strongly F-regular, JT-regular) module, then it is  weakly-endoregular,  weakly-morphic and reduced.    Since all the other forms    are F-regular by  (a), they are weakly-morphic, reduced and  co-reduced as well.
\end{enumerate}
\end{remark}

\begin{corollary}~\label{fagita}
Let $R$ be a commutative ring and $M$ be a nontrivial $R$-module.  Then $M$ is strongly F-regular  whenever  for each $m\in M,M/mR$  is finitely presented and any one  of the following statements is satisfied:
\begin{enumerate}
\item[\emph{(1)}] $mR$ is (weakly-)morphic and reduced,
\item[\emph{(2)}] $mR$ is  Abelian endoregular,
\item[\emph{(3)}] $mR$ is weakly-endoregular,
\item[\emph{(4)}] $mR$ is weakly JT-regular,
\item[\emph{(5)}] $mR$ is co-reduced,
\item[\emph{(6)}] $mR$ is F-regular.
\end{enumerate}
\end{corollary}
\begin{proof}
In view of Proposition~\ref{346} and the fact in Remark~\ref{231x} (b) that the weakly JT-regular modules are the co-reduced modules, it is enough to prove that each one of the given statements (1) to (6)  implies $M$ is F-regular.   Assume that for each  $m\in M,mR$ satisfies any one of the statements given.    Then $M$ is F-regular by Theorem~\ref{230}.
 \end{proof}

\noindent Let $R$ be a commutative ring and $M$ be a nontrivial $R$-module.    Table~\ref{dop} illustrates how the properties:   ``(weakly-)morphic module'' and ``(co-)reduced module''  transfer from a module to each of its cyclic submodules and conversely.   Further, the table shows how  these properties    determine the nature of regularity possessed by a module.

\begin{table}[h!]
  \caption{Regular, (weakly-)morphic and reduced modules and cyclic submodules}~\label{dop}
\[\arraycolsep=1.7pt\def\arraystretch{1.0}
\begin{array}
[c]{|ccccccccc|}\hline
 &                          &                     &                         &                & & && \\
 M~\text{is}&    & \Rightarrow      &   & M~\text{is}&    \Rightarrow      &M~\text{is} &\Rightarrow  &M~\text{is} \\
\text{Strongly}   &    &   & &\text{F-regular}    &           &\text{weakly-endoregular}& &\text{weakly}\\
\text{F-regular}  &                    &&               &               &                &  && \text{JT-regular}\\
 &                          &                     &                         &                & & && \\
\Uparrow &                          &                    &                 & \Updownarrow      &   &\Updownarrow&& ||\\
 &                          &                     &                         &                                                         & &\\
\forall m\in M,&    & \Rightarrow      &   & \forall m\in M,&    \Rightarrow      &M~\text{is} &\Rightarrow   &M~\text{is}  \\
mR~\text{is weakly-morphic} &            & &      &mR~\text{is weakly-morphic} &  &\text{weakly-morphic}     &       &\text{co-reduced}\\
+~mR~\text{is reduced}  &    &    & & +~mR~\text{is reduced}   &       &+~\text{reduced}         &             &    \\
+~M/mR~\text{is~f.p} &         &                       &       &  &                &                      &                 &\\
 &          &                     &                         &                        &                   &\\
\Updownarrow &               &                 &                     &   \Updownarrow&                & \Updownarrow &&||\\
&                          &                       &                     &  &                &                      &                 &\\
\forall m\in M, &    &  \Rightarrow &   &  \forall m\in M,&    \Rightarrow      &M~\text{is} &\Rightarrow &M~\text{is} \\
mR~\text{is co-reduced} &    &                     &   &mR~\text{is co-reduced} &                         &\text{weakly-morphic} &  &\text{co-reduced}  \\
+~M/mR~\text{is~f.p} &    &                     &   &   &                         &+~\text{co-reduced}       &                     &  \\
  &    &                     &   &    &                         &       &                     &  \\\hline
\end{array}
\]
\end{table}

\begin{Ex} The implications in the rows, in general,  cannot be reversed.
\begin{enumerate}[(a)]
\item Co-reduced $\not\Rightarrow$ weakly-morphic.    Let $p$ be a prime element of $\Bbb Z$.   Then the Pr$\ddot{\text{u}}$fer $p$-group $\Bbb Z_{p^\infty}$ is a co-reduced $\Bbb Z$-module.  However, since any non-zero endomorphism of the type $\varphi_a$ of $\Bbb Z_{p^\infty}$ is surjective but not injective,  $\Bbb Z_{p^\infty}$  is  not weakly-morphic as a $\Bbb Z$-module  (see \cite[Proposition 4 and Example 2.2]{kimuli2021}).
\item Weakly-morphic + (co-)reduced on $M \not\Rightarrow$ weakly-morphic  on  cyclic submodules module of $M$.   The $\Bbb Z$-module $\Bbb Q$ is weakly-morphic, (co-)reduced.   However,  its cyclic $\Bbb Z$-submodule $\Bbb Z$ is not weakly-morphic.
\item Co-reduced (= Weakly-morphic + reduced)  on $mR,m\in M \not\Rightarrow M/mR$ is finitely presented.
Teply in \cite{teply1971note}  constructed  a commutative regular ring $R$  with a finitely generated F-regular $R$-module $M$ having
a  submodule $T(M):=\{m\in M :r_R(m)~\text{is an essential ideal of}~R\}.$    By Theorem~\ref{230}, each $mR,m\in M$ is weakly-morphic, reduced and co-reduced.   However, since by \cite{teply1971note}  $T(M)$  is a cyclic pure submodule which is not a direct summand of $M$,   $M/T(M)$ is not finitely presented by Lemma~\ref{345}.
\end{enumerate}
\end{Ex}

\begin{corollary}~\label{trash}
Let $R$ be a commutative ring and $M$ an $R$-module.
$M$ is weakly-endoregular if and only if it is weakly JT-regular and weakly-morphic if and only if it is weakly-morphic and reduced.
\end{corollary}
\begin{proof}
 Since  weakly JT-regular modules are exactly the co-reduced modules, the proof  follows by Theorem~\ref{35}.
\end{proof}

\section{Coincidence of morphic, reduced and regular modules}

\noindent This section  gives conditions under which  the different regularity notions of modules coincide with weakly-morphic and reduced modules.   Further, under some special conditions, we give the kind of regularity a module will attain whenever every (cyclic) submodule of such a module is (weakly-)morphic and reduced.
 Note that (using Lemma~\ref{trey}, Example~\ref{tgvr}, Remark~\ref{231x} and \cite[pg. 15 \& Example 35 (6)]{anderson2019module}) F-regular $\Rightarrow$ weakly-endoregular  $\not\Rightarrow$ F-regular $\Rightarrow$  weakly JT-regular $\not\Rightarrow$  F-regular.

\begin{theorem}~\label{344}
Let $R$ be a commutative ring and $M$ be a nontrivial finitely generated $R$-module.   Then the following statements are equivalent:
\begin{enumerate}
\item[\emph{(1)}]$M$  is weakly-morphic and reduced,
\item[\emph{(2)}]$R/\text{Ann}_R(M)$ is a regular ring,
\item[\emph{(3)}]$M$ is weakly-endoregular,
\item[\emph{(4)}]$M$ is weakly JT-regular,
\item[\emph{(5)}] $M$ is  F-regular,
\item[\emph{(6)}] Every  cyclic submodule  of $M$ is a (weakly-)morphic and reduced (resp., weakly-endoregular, Abelian endoregular,  co-reduced,    weakly JT-regular, F-regular) module.
\end{enumerate}
\end{theorem}
\begin{proof}~
\begin{itemize}
\item[(1)$\Leftrightarrow$(3)] Follows from   Theorem~\ref{35}.
\item[(1)$\Leftrightarrow$(2)] $\Leftrightarrow(4)$ Follows Corollary~\ref{342} and Remark~\ref{231x} (b) respectively.
\item[(3)$\Leftrightarrow$(5)] Follows from \cite[Theorem 22]{anderson2019module}.
\item[(5)$\Leftrightarrow$(6)] Follows from Theorem~\ref{230} and Corollary~\ref{fagita}.
\end{itemize}
 \end{proof}

\noindent Note that the $\Bbb Z$-module $\Bbb Q$ is a non-finitely generated $\Bbb Z$-module that satisfies (1), (3) and (4)  of {\rm Theorem~\ref{344}} but fails on (2), (5) and (6).   Like for rings, the notions of (weakly-)morphic and reduced modules  connect well  to provide  conditions related to  regularity in modules.    In the subcategory of finitely generated modules,  the two properties combined    coincide with different regularity notions in  Theorem~\ref{344}.      Now we give a condition in Proposition~\ref{348} when the endoregular and the strongly F-regular modules coincide with the modules in Theorem~\ref{344}.   Further, we characterise the endoregular and the strongly F-regular modules in terms of  (weakly-)morphic and reduced (sub)modules.

\begin{proposition}~\label{348}
Let $R$ be a commutative ring and $M$ be a nontrivial finitely presented $R$-module.  Then the following statements are equivalent:
\begin{enumerate}
\item[\emph{(1)}]$M$  is weakly-morphic and reduced,
\item[\emph{(2)}]$R/\text{Ann}_R(M)$ is a regular ring,
\item[\emph{(3)}]$M$ is (weakly-)endoregular,
\item[\emph{(4)}]$M$ is weakly JT-regular,
\item[\emph{(5)}]$M$ is  (strongly) F-regular,
\item[\emph{(6)}] Every  cyclic submodule  of $M$ is a (weakly-)morphic and reduced (resp., (weakly-)endoregular, Abelian endoregular,  co-reduced,    weakly JT-regular, F-regular) module.
\end{enumerate}
\end{proposition}
\begin{proof}
The equivalence of $(2)\Leftrightarrow(3)\Leftrightarrow(5)$ follows from \cite[Theorem 23]{anderson2019module}.   The rest of the equivalences follow from  Theorem~\ref{344}.
\end{proof}

\begin{remark}~\label{449}
None of the following notions:  $M$ is ``reduced'', ``weakly-morphic + reduced'', ``weakly-morphic + co-reduced''  implies $S:=\text{End}_R(M)$ is a reduced ring.
Hence, neither weakly JT-regular, (strongly) F-regular, weakly-endoregular implies Abelian endoregular.
There exists a reduced module $M$ with every cyclic submodule weakly-morphic, reduced and co-reduced but with $S$ not reduced, see {\rm Example~\ref{reto}}.
\end{remark}

\begin{Ex}{\rm\cite[Example 24]{anderson2019module}}~\label{reto}
 Let $R$ be a commutative  regular ring with a non-finitely generated maximal ideal $\mathcal{M}$, and let $\overline{R}:=R/\mathcal{M}$ and $M:=R\bigoplus \overline{R}$.  Then $M$ is a finitely generated strongly F-regular module and therefore,  by {\rm Lemma~\ref{trey}}, $M$ is weakly-morphic,  reduced and co-reduced.   However, we claim that $S:=\text{End}_R(M)$ is not a reduced ring.   Note that since
$$S\cong
\left[
  \begin{array}{cc}
    \text{End}_R(R)  & \text{Hom}_R(\bar{R},R) \\
    \text{Hom}_R(R,\bar{R}) & \text{End}_R(\bar{R}) \\
  \end{array}
\right]\cong
 \left[
                                                                                                    \begin{array}{cc}
                                                                                                      R & 0 \\
                                                                                                      \overline{R} & \overline{R} \\
                                                                                                    \end{array}
                                                                                                  \right],
    $$
the  endomorphism $\varphi$ corresponding to $\left[
                                                         \begin{array}{cc}
                                                            0 & 0 \\
                                                            \overline{1} & \overline{0} \\
                                                          \end{array}
                                                        \right]
    $ is non-zero but  $\varphi^2=0$.

\end{Ex}

\noindent
\begin{definition}~\label{trees}
Let $R$ be a commutative ring.  An  $R$-module  $M$ is   {\it almost locally simple module} \cite{anderson2019module}  if $M_{\mathcal{M}}$ is a trivial or  simple $R_{\mathcal{M}}$-module (equivalently if $M_{\mathcal{M}}$ is  a trivial or simple $R$-module) for each maximal ideal $\mathcal{M}$ of $R$.
\end{definition}

\noindent It is well known that $R$ is  an almost locally simple $R$-module if and only if $R$ is a regular ring.  By Anderson, Chun \& Juett in \cite[pg. 2]{anderson2019module}, the  ``almost locally simple property'' in modules is another form of module-theoretic regularity.

\begin{lemma}~\label{sothat} Let $R$ be a commutative ring and $M$ be a nontrivial $R$-module.  Then
\begin{enumerate}
\item[\emph{(1)}] {\rm\cite[Theorem 4]{anderson2019module}}~\label{sothau}
$M$ is JT-regular if and only if $M$ is a multiplication and weakly JT-regular module;
\item[\emph{(2)}]{\rm\cite[Theorem 13]{anderson2019module}} $M$ is JT-regular $ \Rightarrow M$ is almost locally simple $ \Rightarrow M$ is strongly F-regular $ \Rightarrow M$ is F-regular $ \Rightarrow M$ is weakly JT-regular.
\end{enumerate}
\end{lemma}

\begin{corollary}~\label{titles}
Let $R$ be a commutative ring and $M$ be a nontrivial multiplication $R$-module.  Then
$M$ is JT-regular $ \Leftrightarrow M$ is almost locally simple $ \Leftrightarrow M$ is strongly F-regular $ \Leftrightarrow M$ is F-regular $ \Leftrightarrow M$ is weakly JT-regular.
\end{corollary}

\begin{proof}
This is immediate from Lemma~\ref{sothat}.
\end{proof}
\noindent Naoum in \cite{naoum1995regular} proved that a multiplication module is  strongly F-regular if and only if its ring of endomorphisms $S$ is regular (i.e., $M$ is endoregular).
Proposition~\ref{349} shows that for finitely generated multiplication modules over commutative rings, the (weakly-)morphic and reduced modules coincide with all the regularity notions we have discussed.

\begin{proposition}~\label{349}
Let $R$ be a commutative ring and $M$ be a nontrivial finitely generated multiplication $R$-module.  Then the following statements are equivalent:
\begin{enumerate}
\item[\emph{(1)}]$M$ is (weakly-)morphic and reduced,
\item[\emph{(2)}]$R/\text{Ann}(M)$ is  a regular ring,
\item[\emph{(3)}]$M$ is (weakly-)endoregular,
\item[\emph{(4)}]$M$ is (Abelian) endoregular,
\item[\emph{(5)}]$M$ is (weakly) JT-regular,
\item[\emph{(6)}]$M$ is almost locally simple,
\item[\emph{(7)}]$M$ is  (strongly) F-regular,
\item[\emph{(8)}] Every  cyclic submodule  of $M$ is a (weakly-)morphic and reduced (resp., (weakly-)endoregular, Abelian endoregular,  co-reduced,    (weakly) JT-regular, (strongly) F-regular,  almost locally simple) module.
\end{enumerate}
\end{proposition}
\begin{proof}~
\begin{itemize}
\item[(1)$\Leftrightarrow$(2)] $\Leftrightarrow$ (3) $\Leftrightarrow$ (5) $\Leftrightarrow$ (8) Since  weakly-morphic finitely generated multiplication modules are morphic by Lemma~\ref{300}, the proof of the equivalence follows from
 Theorem~\ref{344} and Corollary~\ref{titles}.
 \item[(5)$\Leftrightarrow$(6)] $\Leftrightarrow$ (7) Follows from Corollary~\ref{titles}.
\item[(3)$\Leftrightarrow$(4)] Follows from  Proposition~\ref{444}.
\end{itemize}
 \end{proof}

\noindent
 Ware   \cite[Definition 2.3]{ware1971endomorphism} calls  $M$ {\it regular} (call it {\it W-regular}) if it is projective and every homomorphic image of $M$ is flat, or equivalently if $M$ is projective and every cyclic submodule of $M$ is a direct summand.
To extend the Ware regularity notion,  Zelmanowitz   defined the non-projective regular modules.   $M$ is  a {\it Zelmanowitz regular} \cite{zelmanowitz1972regular} (call it {\it Z-regular}) module  if given any $m\in M$, there exists $\varphi\in\text{Hom}_R(M,R)$ such that $m\varphi(m) = m$, or equivalently if for any $m\in M, mR$ is projective and is a direct summand of $M$.

\begin{remark} Let $R$ be a commutative ring and $M$ be a nontrivial $R$-module.
\begin{enumerate}[(a)]

\item Since (by \cite{ramamurthi1973finitely}) $M$ is W-regular $ \Rightarrow M$ is Z-regular $ \Rightarrow M$ is  strongly F-regular $\Rightarrow M$ is F-regular, the W-regular modules and the Z-regular modules are  weakly-endoregular, weakly-morphic and reduced by Lemma~\ref{trey}.
\item If $M$ is projective, then $M$ is W-regular $ \Leftrightarrow M$ is Z-regular $ \Leftrightarrow M$ is  (strongly) F-regular $\Leftrightarrow$ for each $m\in M, mR$ is a (weakly-)morphic and reduced module $\Leftrightarrow$ for each $m\in M, mR$ is a co-reduced module (\cite[Proposition 2]{ware1971endomorphism} and Theorem~\ref{230}).
\end{enumerate}
\end{remark}

\begin{ack} We would like to thank the anonymous referee for a careful proofreading of our manuscript and  for making many valuable comments.   This work was carried out at Makerere University with support from ISP through the Makerere-Sida Bilateral Program Phase IV, Project 316 ``Capacity Building in Mathematics and its Application" and through the Eastern Africa Algebra Research Group (EAALG).
\end{ack}

\end{document}